\documentclass{article}
\usepackage{spconf,amsmath,graphicx,comment, amsfonts, url}

\RequirePackage[colorlinks,citecolor=blue,urlcolor=blue,linkcolor=blue]{hyperref}

\DeclareMathOperator*{\argmin}{\arg\!\min}

\title{3D {\em ab initio} modeling in cryo-EM by autocorrelation analysis}
\name{Eitan Levin, Tamir Bendory, Nicolas Boumal, Joe Kileel and Amit Singer\thanks{This work was supported by Award Number R01GM090200 from the NIGMS, the Simons Foundation Investigator Award, the Simons Collaboration in Algorithms and Geometry, the Moore Foundation Data-Driven Discovery Investigator Award and NSF DMS-1719558.}}
\address{Princeton University, Program in Applied and Computational Mathematics, Princeton NJ, USA
}

\begin{document}
%
\maketitle
\begin{abstract}
Single-Particle Reconstruction (SPR) in Cryo-Electron Microscopy (cryo-EM) is the task of estimating the 3D structure of a molecule from a set of noisy 2D projections, taken from unknown viewing directions. Many algorithms for SPR start from an initial reference molecule, and alternate between refining the estimated viewing angles given the molecule, and refining the molecule given the viewing angles. This scheme is called iterative refinement.  Reliance on an initial, user-chosen reference introduces model bias, and poor initialization can lead to slow convergence. Furthermore, since no ground truth is available for an unsolved molecule, it is difficult to validate the obtained results.  This creates the need for high quality \emph{ab initio} models that can be quickly obtained from experimental data with minimal priors, and which can also be used for validation. We propose a procedure to obtain such an \emph{ab initio} model directly from raw data using Kam's autocorrelation method.  Kam's method has been known since 1980, but it leads to an underdetermined system, with missing orthogonal matrices. Until now, this system has been solved only for special cases, such as highly symmetric molecules or molecules for which a homologous structure was already available.  In this paper, we show that knowledge of just two clean projections is sufficient to guarantee a unique solution to the system. This system is solved by an optimization-based heuristic. For the first time, we are then able to obtain a low-resolution \emph{ab initio} model of an asymmetric molecule directly from raw data, without 2D class averaging and without tilting. Numerical results are presented on both synthetic and experimental data.
\end{abstract}
\section{Introduction}
\label{sec:intro}
Cryo-EM is an increasingly popular method for determining the 3D structure of molecules, especially those that resist crystalization~\cite{resolution_revolution, Revol, Frank2006}. Advances in this technique were recognized by the 2017 Nobel Prize in Chemistry \cite{Nobel2017}. For SPR in cryo-EM, a sample containing many (ideally) identical molecules in unknown orientations are frozen in a sheet of ice. An electron microscope produces a top view of the sample in one image, called a micrograph, from which projection images of individual molecules are extracted in a process called particle picking. In order to limit radiation damage to the organic molecules caused by the electron beam, the electron dosage must be kept low, resulting in a low signal-to-noise ratio (SNR) in each of the projections. In addition, the images are affected by the Contrast Transfer Function (CTF) of the microscope, causing further aberrations. The goal is to estimate the 3D structure of the molecule from a large set of projections selected from multiple micrographs.

Typical approaches to SPR use iterative refinement procedures that start from an initial guess of the 3D structure, apply a low-pass filter, and then refine it by alternating between estimation of the viewing directions of the projections given the molecule and vice versa~\cite{relion, Barnett2016, Punjani2017}. Since these algorithms solve a non-convex problem, the quality of their output as well as the speed of their convergence depend on the initialization, particularly at low SNR or with small particles \cite{Bhamre2015, Bhamre2017}.

In contrast, \emph{ab initio} methods do not require an initial model. Currently, few \emph{ab initio} methods are available. The random conical tilt method~\cite{Radermacher2} requires the molecule to have a strongly preferred orientation. Methods that do not involve tilting are either based on moments~\cite{Salzman1990, Goncharov1988} or common lines~\cite{Vainshtein, VanHeel}. However, these approaches typically fail to recover the 3D structure from non-averaged experimental images due to the low SNR.

We present a new method called {\em orthogonal matrix retrieval by projection matching}, based on Kam's autocorrelation analysis \cite{kam1980,Kam1985}. Unlike the above mentioned methods for \emph{ab initio} modeling, Kam's method completely sidesteps estimation of particle orientations. It only requires the covariance matrix of the projection images, which can be estimated accurately for any SNR given sufficiently many particle images. Kam's analysis recovers the expansion coefficients of the structure, up to a sequence of missing orthogonal matrices. It assumes the viewing directions are uniformly distributed over the sphere. Recently, there have been numerous attempts to apply Kam's method to XFEL \cite{Saldin2011, Saldin2010, Kurta2017} and to cryo-EM \cite{Bhamre2015, Bhamre2017}. Restrictingly, the first make either strong symmetry assumptions on the molecule or limit the rotations to a single axis, while the latter assume that the structure of a similar molecule is already available.

In this work, we apply Kam's method to resolve the molecular structure directly from raw experimental images without estimating viewing directions, for the first time. We use the method of \cite{Bhamre2016} to estimate the covariance matrix of the projections from raw data. We then recover the missing orthogonal matrices by matching to two clean or denoised images, via Riemannian optimization.
The computational complexity of our algorithm is linear in the number of images.
As an information-theoretic guarantee, we prove that 2D covariance together with merely two clean images uniquely determine the 3D molecular structure.
For reproducibility, a Matlab implementation of our method is available at \url{https://github.com/eitangl/kam_cryo}.

The rest of this paper is organized as follows. In Section~\ref{sec:img_form_model}, we describe the image formation model in cryo-EM. In Section~\ref{sec:general_pipeline}, we describe Kam's autocorrelation analysis and formulate the orthogonal matrix retrieval problem. Section~\ref{sec:orth_mat_ret_problem} describes our procedure for solving the orthogonal matrix retrieval problem, which enables us to recover the molecular structure, and provides an information-theoretic guarantee.  In Section~\ref{sec:numer_examples}, we show the performance of our method on synthetic and experimental datasets. Finally, in Section~\ref{sec:limits}, we discuss possible extensions of the method for future work.

\section{Image formation model}\label{sec:img_form_model}
Let $\phi : \mathbb{R}^{3} \rightarrow \mathbb{R}$ be the Coulomb potential representing the molecular structure we wish to estimate. The $j^{\textup{th}}$ projection image $I_j : \mathbb{R}^{2} \rightarrow \mathbb{R}$ is modeled as
\begin{equation}\label{eq:img_formation}
I_j = H_j\ast\mathcal{P}_j[\phi] \,+\, \varepsilon_j, \quad j = 1, \ldots, n.
\end{equation}
Here, $H_j \! : \mathbb{R}^{2} \rightarrow \mathbb{R}$ corresponds to the CTF affecting the $j^{\textup{th}}$ image by convolution, $\varepsilon_j$ is noise and $\mathcal{P}_j$ is the tomographic projection operator given by
\begin{equation}\mathcal{P}_j[\phi](x,y) = \int_{-\infty}^{\infty}\phi(R_j^Tr) \,\, dz,\end{equation}
where $r=(x,y,z)^T$ are Cartesian coordinates and $R_j \in \textup{SO}(3)$ is the orientation of the $j^{\textup{th}}$ particle. This formation model is more neatly expressed in the Fourier domain. Owing to the Fourier-slice theorem~\cite[pp. 11]{Natterer2001a}, the 2D Fourier transform of a 2D projection image is the restriction of the 3D Fourier transform of $\phi$ to the plane passing through the origin perpendicular to the viewing direction. Denoting Fourier transforms by hats, we can rewrite the formation model as
\begin{equation}\label{eq:img_form_Fourier} \widehat I_j(k_x, k_y) = \widehat H_j\cdot\widehat\phi\big{(} R_j^T (k_x, k_y, 0)^T \big{)}  + \, \widehat\varepsilon_j,\end{equation}
where $k_x, k_y$ are Cartesian coordinates in 2D Fourier space.

\section{Kam's autocorrelation analysis}\label{sec:general_pipeline}
We assume that the structure $\phi$ is essentially compactly supported and bandlimited with bandlimit $c$.  We expand the Fourier transform of the density in the eigenfunctions of the Laplacian with Dirichlet boundary conditions over the radius $c$ ball in $\mathbb{R}^3$, working in spherical coordinates:
\begin{equation}\label{eq:vol_expansion}
\widehat\phi(k,\theta,\varphi) = \sum_{l=0}^{L}\sum_{m=-l}^l\sum_{s=1}^{S(l)}a_{lms}Y_{lm}(\theta,\varphi)j_{ls}(k). \end{equation}
Here, $Y_{lm}$ are the real spherical harmonics and
\begin{equation}\label{eq:radial_funcs} j_{ls}(k) = \frac{\sqrt{2}}{c^{3/2}|j_{l+1}(u_{l,s})|}\, j_l(u_{l,s}k/c),\end{equation}
where $j_l$ is the spherical Bessel function of order $l$, the scalar $u_{l,s}$ is the $s^{\textup{th}}$ positive zero of $j_l$.

We shall assume a bandlimit $c$ smaller than the Nyquist frequency
and a finite expansion of the above form, since we focus on recovering a low-resolution version of the molecule, suitable for an \textit{ab initio} estimate. The truncation limit $S(l)$ is chosen by the sampling criterion proposed in~\cite[Eq.~(8)]{Bhamre2017}, enforcing essentially compact support in real space, and $S(l)$ is a monotonically decreasing function of $l$.

Our goal is to estimate the coefficients $a_{lms}$. In a seminal paper~\cite{kam1980}, Kam showed that the matrices
\begin{equation}\label{eq:Cl_mats} C_l(s_1,s_2) = \sum_{m=-l}^la_{lms_1}\overline{a_{lms_2}},\quad l=0,\ldots, L,\end{equation}
can be recovered directly from the noisy projections, provided that the viewing directions are uniformly distributed over the sphere.

Defining the $S(l)\times (2l+1)$ matrix of coefficients $A_l$ indexed as $A_l(s,m)=a_{lms}$ for fixed $l$, the $S(l)\times S(l)$ matrices $C_l$ in Eq.~(\ref{eq:Cl_mats}) satisfy the relation
\begin{equation}\label{eq:Cl_chol} C_l = A_l^{}A_l^*, \end{equation}
where $A^*$ denotes the Hermitian conjugate of $A$.
Since the molecular density $\phi$ is real-valued, its Fourier transform is conjugate-symmetric, and hence the matrices $A_l$ are purely real for even $l$, and purely imaginary for odd $l$. Therefore, Eq.~(\ref{eq:Cl_chol}) determines $A_l$ uniquely up to an orthogonal matrix of size $(2l+1)\times(2l+1)$.

Formally, we take a Cholesky decomposition of the estimated $C_l$ to obtain $S(l)\times(2l+1)$ matrices $F_l$. Accordingly, $A_l=F_lO_l$ for some unknown $(2l+1)\times(2l+1)$ orthogonal matrices $O_l$. This is the missing orthogonal matrix problem in Kam's method, which we aim to solve with minimal priors on the molecule. This would then allow us to recover the 3D structure.

\section{Orthogonal matrix retrieval by projection matching}\label{sec:orth_mat_ret_problem}

We begin by noting that the matrix $O_0$ is just a sign $\pm1$, and can be easily recovered from the average image of the dataset. Specifically, we take the radially-isotropic average of all the projections, and note that this average is determined only by the $l=0$ component, proportional to $\sum_{s=1}^{S(0)}a_{00s}j_{0s}(k)$. This determines the coefficients for $l=0$.

The main contribution of this paper is the observation that the remaining $\{O_l\}_{l=1}^L$ may be retrieved by matching to merely two clean or denoised projections. These projections can be obtained for example by using the Wiener filter-based method of \cite{Bhamre2016} to denoise and CTF-correct individual projection images. To see how a known clean projection constrains the missing $\{O_l\}$, we write Eq.~(\ref{eq:vol_expansion}) in matrix form to get
\begin{equation}\label{eq:linear_dep_orth_mats}
\widehat\phi(\{O_l\}) = \sum_{l=0}^Lj_l F_lO_l Y_l,
\end{equation}
where matrices are indexed as $[j_l]_{k, s} = j_{ls}(k)$, $[Y_l]_{m,(\theta,\varphi)}=Y_{lm}(\theta,\varphi)$ and $F_l$ is obtained from a Cholesky decomposition of $C_l$ mentioned in Section~\ref{sec:general_pipeline}.

Without loss of generality, the first clean projection is the restriction to the $k_{x}k_{y}$-plane of $\widehat\phi$, as the molecule can only be estimated up to a global rotation and reflection. Now, let the orientation of the particle in the second clean image be given by an unknown $R \in \textup{SO}(3)$.  Since rotating real spherical harmonics of degree $l$ may be expressed with matrix multiplication \cite{Ivanic1996}, the second clean image also imposes linear constraints on $\{O_l\}$.  Writing $D_l^{(R)}$ for the $(2l + 1) \times (2l + 1)$ Wigner D-matrix of $R$ in this irreducible representation of $\textup{SO}(3)$, the second projection imposes the constraints
\begin{equation}\label{eq:rotated_plane_constraint}
\widehat\phi(\{O_l\}, R) = \sum_{l=0}^L j_l F_l O_l D_l^{(R)} Y_l.
\end{equation}

From Rodrigues' formula for associated Legendre polynomials, restricting $Y_{lm}$ to $\theta=\pi/2$ (the $k_{x}k_{y}$-plane) sets all the rows of $Y_l$ for which $l\not\equiv m \,\, (\textup{mod } 2)$ to zero.  Thus the clean projections constrain every other column of $\{O_l\}$ and $\{O_l D_l^{(R)}\}$, respectively. It can be shown (proof omitted here due to space limitations) that under mild technical conditions these linear constraints in fact uniquely determine the orthogonal matrices $\{O_l\}$ and the rotation $R$, and hence the 3D structure itself.

In practice, given two clean or denoised images $I^{c}_{1}, I^{c}_{2}$, we begin by matching to each image separately.  To do this, we assume both projections lie on the $k_xk_y$-plane corresponding to the $3\times 3$ identity rotation matrix $\mathcal{I}_{3}$, let $\widehat\phi(\{O_l\})_{k_xk_y}$ denote the restriction of Eq.~\ref{eq:linear_dep_orth_mats} to the $k_xk_y$ plane, and obtain estimates for every other column in two sets of orthogonal matrices:
\begin{equation}\label{eq:orth_mat_optim_prob} \begin{aligned}
\{o_{l;1}\}_{l=1}^{L} &= \argmin_{\substack{o_{l} \in \mathbb{R}^{(2l+1) \times (l+1)}\\o_{l}^{T} o_{l}^{} = \mathcal{I}_{l+1}}} ||\widehat\phi\big{(}\{o_{l}\}\big{)}_{k_xk_y} - \widehat{I^{c}_1}||_F^2,
\\
\{o_{l;2}\}_{l=1}^{L} &= \argmin_{\substack{o_{l} \in \mathbb{R}^{(2l+1) \times (l+1)}\\o_{l}^{T} o_{l}^{} = \mathcal{I}_{l+1}}} ||\widehat\phi\big{(}\{o_{l}\}\big{)}_{k_xk_y} - \widehat{I^{c}_2}||_F^2.
\end{aligned}\end{equation}
We estimate $\{o_{l;1}\}, \{o_{l;2}\}$ via optimization over the appropriate product of manifolds using {Manopt} \cite{manopt}.
Note that Riemannian gradient descent is only guaranteed to converge to critical points of the cost function~\cite{AMS08}. However, for the purposes of \textit{ab initio} modeling, our implementation performs satisfactorily, as seen empirically in Section \ref{sec:numer_examples}.

Continuing the algorithm, we then merge results from the two images together.  Writing $O_{l}$ for the missing orthogonal matrices, taking the first image to have identity orientation and the second $R$, it follows that every other column of $O_{l}$ should equal columns of $o_{l;1}$ while every other column of $O_{l}D_l^{(R)}$ should equal columns of $o_{l;2}$.  We solve for $R$ and $\{O_{l}\}$ by making these as consistent as possible.
Formally, for each $l=1,\ldots, L$, we form the matrices $D_l = \left[\widetilde{\mathcal{I}}_{2l+1}\, |\, \widetilde D_l^{(R)}\right]$ and $B_l = \left[o_{l; 1} \, | \, o_{l;2} \right]$, where we denote by $\widetilde A\in\mathbb{R}^{(2l+1)\times(l+1)}$ the matrix obtained from $A\in\mathbb{R}^{(2l+1)\times(2l+1)}$ by taking only every other column including the first and last, and where $[X \, | \, Y]$ denotes the horizontal concatenation of matrices $X$ and $Y$. We then solve the minimization
\begin{equation}\label{eq:init_prob}\begin{aligned}
\min_{\substack{R \in \textup{SO}(3)}} \,\, \sum_{l = 1}^{L} \,\, \min_{\substack{{O_{l} \in O(2l+1)}}} \,\, || O_{l}D_{l} - B_{l} ||^{2}_{F}.
\end{aligned}\end{equation}
This is done by densely sampling $R \in \textup{SO}(3)$ and noting that for $R$ fixed the minimization $\min_{\substack{{O_{l} \in O(2l+1)}}} \,\, || O_{l}D_{l} - B_{l} ||^{2}_{F}$ is an instance of the orthogonal Procrustes problem, which has a closed form solution via SVD of $B_{l}^{}D_{l}^{T}$~\cite{Schonemann1966}. Finally, we further refine our estimates of $\{O_l\}$ and $R$ using Manopt.

\section{Numerical Examples}\label{sec:numer_examples}
\begin{figure}[htb]
\begin{minipage}[b]{.6\linewidth}
  \centering
  \centerline{\includegraphics[width=5.5cm]{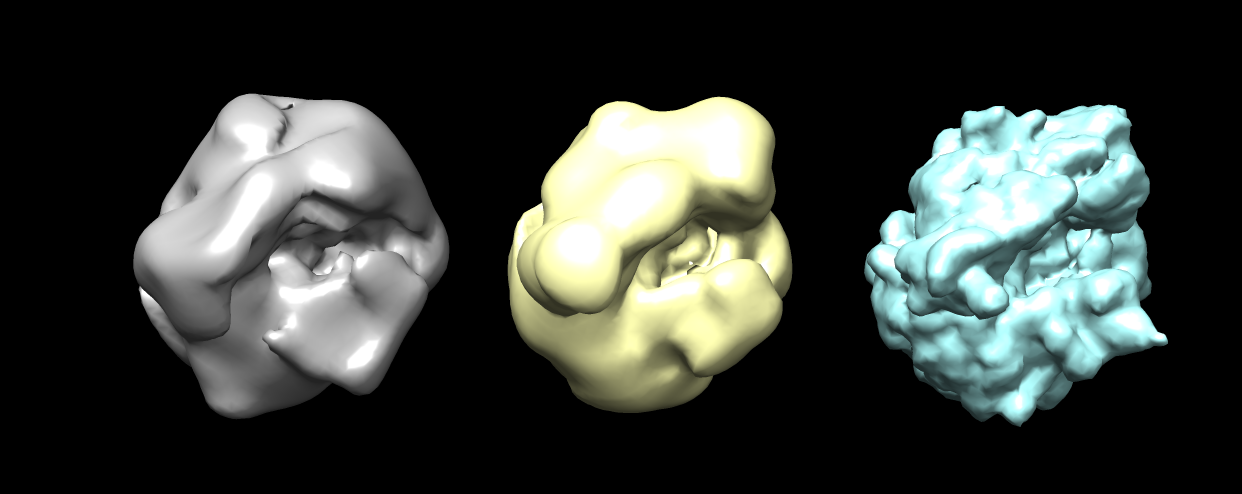}}
  \centerline{(a)}\medskip
\end{minipage}
\hfill
\begin{minipage}[b]{.35\linewidth}
  \centering
  \centerline{\includegraphics[width=3.2cm]{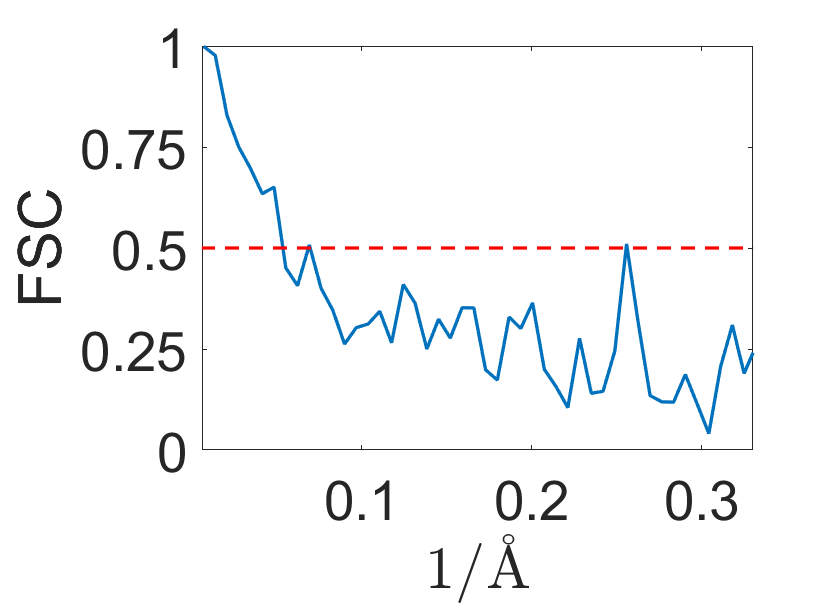}}
  \centerline{(b)}\medskip
\end{minipage}

\begin{minipage}[b]{.6\linewidth}
  \centering
  \centerline{\includegraphics[width=5.5cm]{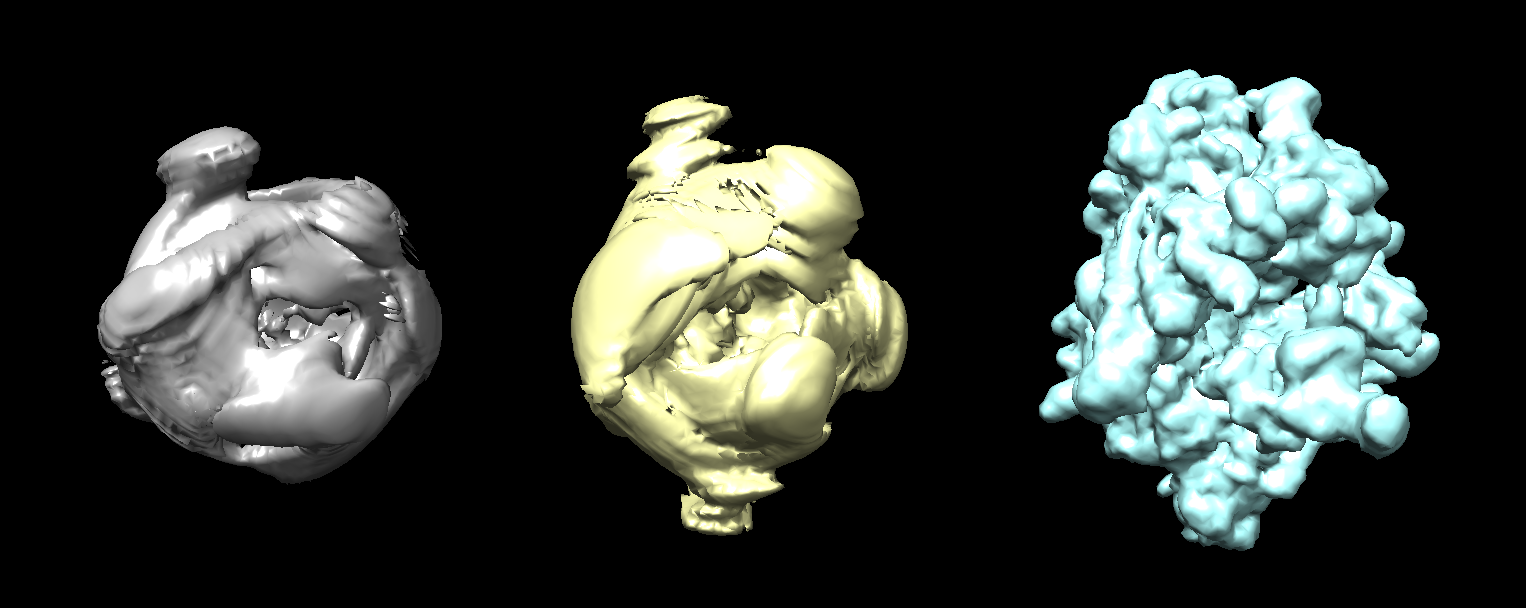}}
  \centerline{(c)}\medskip
\end{minipage}
\hfill
\begin{minipage}[b]{0.35\linewidth}
  \centering
  \centerline{\includegraphics[width=3.2cm]{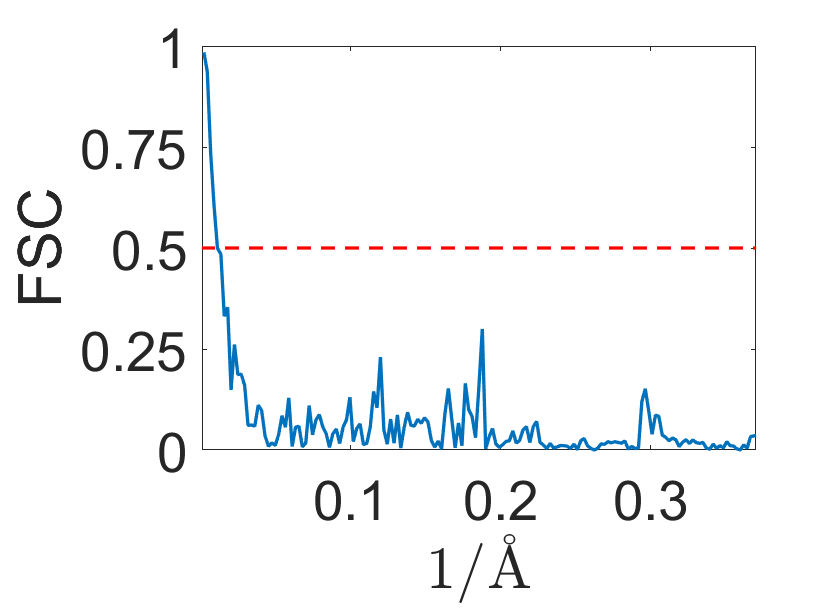}}
  \centerline{(d)}\medskip
\end{minipage}
\caption{\textit{Reconstruction results.} (a) Synthetic data results: the reconstruction (grey), the ground truth (yellow), and the original (blue). (b) FSC curve for synthetic data. (c) Raw data results for yeast mitochondrial
ribosome (EMPIAR-10107): the reconstruction (grey), the ground truth, taken as the corresponding low-resolution EMD-3551 (yellow) and the original EMD-3551, reconstructed using RELION, slightly low-passed with a Gaussian filter to remove noise artifacts (blue). (d) FSC curve for raw data. Note that the reconstruction was measured with respect to the low-resolution ground truth in both (b) and (d).}
\label{fig:res}
\end{figure}

We begin with results on a synthetic dataset consisting of $5\times 10^4$ noisy projections of size $109 \times 109$ with $\text{SNR} = 1/10$ from uniformly random viewing directions of the 70S ribosome with P-site tRNA, available in the Protein Data Bank (EMDB) as EMD-5360. The images are divided into 100 defocus groups, and are centered. The bandlimit is assumed to be $c=1/4$ (half the Nyquist frequency), and the truncation for the expansion of the structure is set to $L=10$. The two images for the reconstruction were chosen uniformly at random. The reconstruction results are presented in Fig.~\ref{fig:res} (a) and (c), where we also show the Fourier Shell Correlation (FSC)~\cite{Harauz1986} of our reconstruction with the low-resolution ground truth. The resolution of the reconstruction is 19 {\AA} using the FSC $=0.5$ criterion.

We also present results on an experimental dataset consisting of $~3.5\times 10^5$ projections of size $330\times 330$ of the yeast mitochondrial ribosome, available in the Electron Microscopy Public Image Archive (EMPIAR) as EMPIAR-10107, out of which we chose $2\times 10^5$ random projections for implementation reasons. We pre-processed the data only by whitening the projections using the method described in \cite[Sec.~2.2]{Bhamre2016}, and directly applied the method of \cite{Bhamre2016} to estimate the covariance matrix of the projections, from which we obtained $\{C_l\}$. Here we set $c=1/4$ and $L=7$. Once again, the two projections for the reconstruction were chosen randomly. We ran the algorithm on a machine with 60
cores, running at 2.3 GHz, with total RAM of 1.5TB. The pre-processing then took $\sim$5 hours, while the reconstruction itself took $\sim$15 minutes.
The resolution of the reconstruction is 89 {\AA} using the FSC $=0.5$ criterion. For the ground truth, we took the RELION reconstruction available as EMD-3551, and expanded it in a finite expansion of the form Eq.~\ref{eq:vol_expansion} with the same truncation as for our own reconstruction. The original EMD-3551 is presented alongside its finite expansion for comparison, slightly low-passed with a Gaussian filter to remove noise artifacts present in the reconstruction.

\section{Conclusions and future work}\label{sec:limits}
In this paper, we have presented a new method to obtain \emph{ab initio} low-resolution 3D molecular structures directly from raw cryo-EM data. We rely on Kam's autocorrelation analysis, which recovers the expansion coefficients of the molecule from the covariance matrix of its projections, up to a set of missing orthogonal matrices. We retrieve these matrices using two clean or denoised projections, and showed that two clean projections determine the structure under mild assumptions. Finally, we demonstrated the performance of the method on both synthetic and experimental datasets. This is the first application of Kam's method to {\em ab initio} modeling of asymmetric molecules from raw experimental data without any averaging.

Nevertheless, we observe in practice that our method is only capable of recovering a low-resolution version of the molecule. While sufficient to initialize iterative refinement algorithms and validate their output, we would like to improve the method to obtain higher-resolution reconstructions, while keeping the computational cost low.
We believe our resolution limitation stems from several features present in real datasets, but which our formulation currently ignores. First, because individual projection images are picked from extremely noisy micrographs, the projections may not be centered. Second, while we assume all the molecules in the sample are identical, in practice they may appear in different conformations, and so the projections would come from several different molecules. Third, Kam's method as stated here assumes the viewing angles of the projections are uniformly distributed over the sphere. In practice, molecules have preferred orientations, which skews the distribution of viewing angles. For symmetric molecules, it can be shown that a single clean image is sufficient to determine the missing orthogonal matrices, in which case our method may be simplified and improved.
We intend to extend Kam's method to account for these features.
Finally, it may be possible to avoid the need for clean projections in the first place by using higher-order correlations in addition to the covariance matrix, as originally suggested by Kam \cite{kam1980}.

\section{Acknowledgements}
We thank Joakim And\'{e}n, Tejal Bhamre, and Yoel Shkolniskly for productive discussions and help with the code.

\begin{figure}[h!]
\begin{minipage}[b]{.22\linewidth}
  \centering
  \centerline{\includegraphics[width=2.5cm]{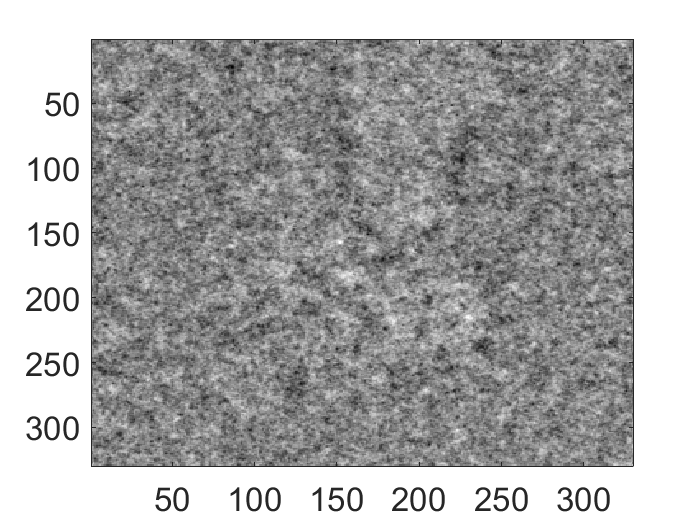}}
  \centerline{(a)}\medskip
\end{minipage}
\hfill
\begin{minipage}[b]{.22\linewidth}
  \centering
  \centerline{\includegraphics[width=2.5cm]{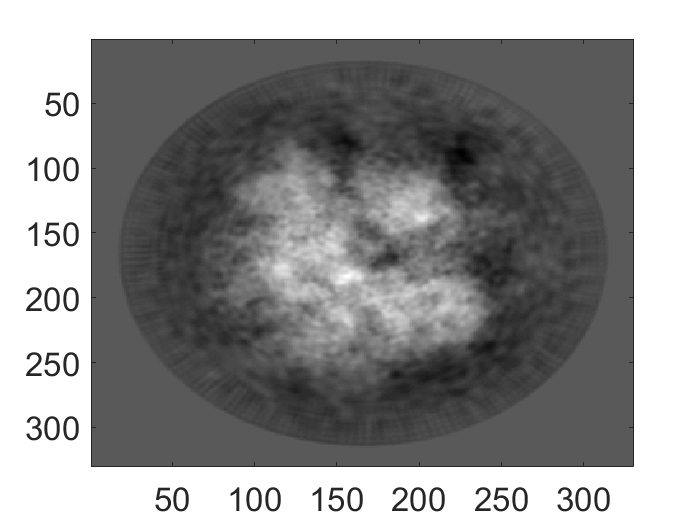}}
  \centerline{(b)}\medskip
\end{minipage}
\hfill
\begin{minipage}[b]{.22\linewidth}
  \centering
  \centerline{\includegraphics[width=2.5cm]{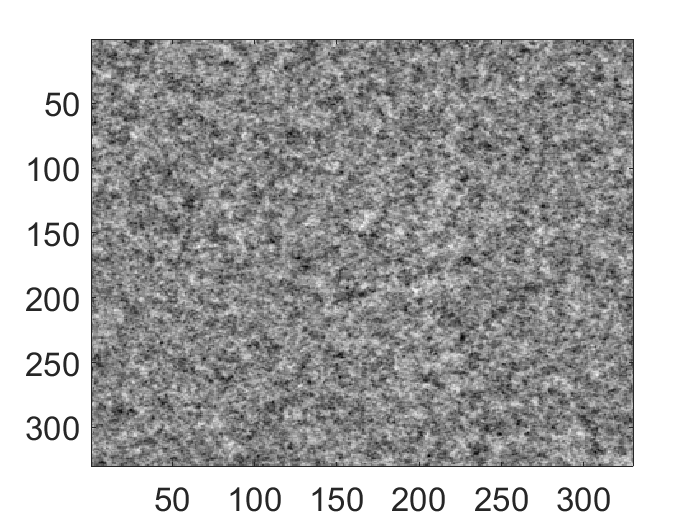}}
  \centerline{(c)}\medskip
\end{minipage}
\hfill
\begin{minipage}[b]{0.22\linewidth}
  \centering
  \centerline{\includegraphics[width=2.5cm]{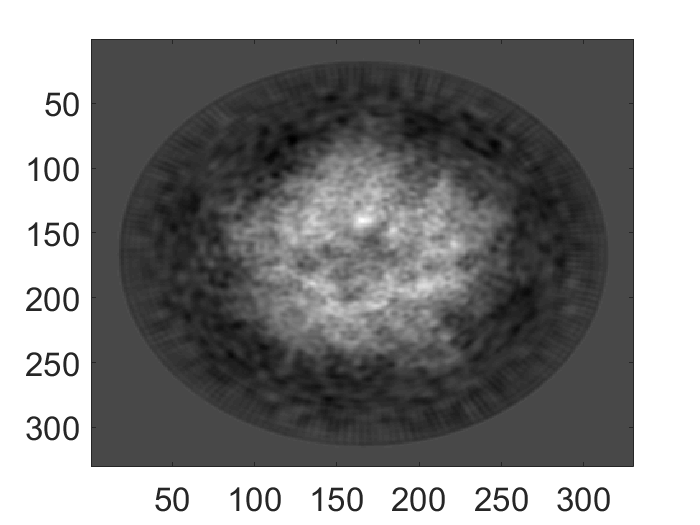}}
  \centerline{(d)}\medskip
\end{minipage}
\caption{\textit{Images used for raw data reconstruction.} (a) and (c) are the original noisy images and (b) and (d) are the corresponding denoised images used for the optimization.}
\label{fig:projs_raw_res}
\end{figure}

\vfill\pagebreak
\bibliographystyle{habbrv}
\bibliography{ISBI_refs}

\begin{thebibliography}{10}

\bibitem{Nobel2017}
\url{https://www.nobelprize.org/nobel_prizes/chemistry/laureates/2017/}.

\bibitem{AMS08}
P.-A. Absil, R.~Mahony, and R.~Sepulchre.
\newblock {\em {Optimization Algorithms on Matrix Manifolds}}.
\newblock Princeton University Press, Princeton, NJ, 2008.

\bibitem{Revol}
X.~Bai, G.~McMullan, and S.~H.~W. Scheres.
\newblock How {cryo-EM} is revolutionizing structural biology.
\newblock {\em Trends Biochem. Sci.}, 40(1):49 -- 57, 2015.

\bibitem{Barnett2016}
A.~Barnett, L.~Greengard, A.~Pataki, and M.~Spivak.
\newblock {Rapid Solution of the Cryo-EM Reconstruction Problem by Frequency
  Marching}.
\newblock {\em SIAM J. Imaging Sci.}, 10(3):1170--1195, 2017.

\bibitem{Bhamre2015}
T.~Bhamre, T.~Zhang, and A.~Singer.
\newblock {Orthogonal matrix retrieval in cryo-electron microscopy}.
\newblock In {\em 2015 IEEE 12th Int. Symp. Biomed. Imaging}, pages 1048--1052.
  IEEE, 2015.

\bibitem{Bhamre2016}
T.~Bhamre, T.~Zhang, and A.~Singer.
\newblock {Denoising and covariance estimation of single particle cryo-EM
  images}.
\newblock {\em J. Struct. Biol.}, 195(1):72--81, 2016.

\bibitem{Bhamre2017}
T.~{Bhamre}, T.~{Zhang}, and A.~{Singer}.
\newblock {Anisotropic twicing for single particle reconstruction using
  autocorrelation analysis}.
\newblock {\em ArXiv e-prints}, 2017, 1704.07969.

\bibitem{manopt}
N.~Boumal, B.~Mishra, P.-A. Absil, and R.~Sepulchre.
\newblock {{M}anopt, a {M}atlab Toolbox for Optimization on Manifolds}.
\newblock {\em J. Mach. Learn. Res.}, 15:1455--1459, 2014.

\bibitem{Frank2006}
J.~Frank.
\newblock {\em {Three-Dimensional Electron Microscopy of Macromolecular
  Assemblies: Visualization of Biological Molecules in Their Native State}}.
\newblock Oxford University Press, New York, 2006.

\bibitem{Goncharov1988}
A.~B. Goncharov.
\newblock Integral geometry and three-dimensional reconstruction of randomly
  oriented identical particles from their electron microphotos.
\newblock {\em Acta Appl. Math.}, 11(3):199--211, 1988.

\bibitem{Harauz1986}
G.~Harauz and M.~van Heel.
\newblock {Exact filters for general geometry three dimensional
  reconstruction}.
\newblock {\em Optik (Stuttg).}, 73:146--156, 1986.

\bibitem{Ivanic1996}
J.~Ivanic and K.~Ruedenberg.
\newblock {Rotation Matrices for Real Spherical Harmonics. Direct Determination
  by Recursion}.
\newblock {\em J. Phys. Chem.}, 100(15):6342--6347, 1996.

\bibitem{kam1980}
Z.~Kam.
\newblock The reconstruction of structure from electron micrographs of randomly
  oriented particles.
\newblock {\em J. Theor. Biol.}, 82(1):15 -- 39, 1980.

\bibitem{Kam1985}
Z.~Kam and I.~Gafni.
\newblock Three-dimensional reconstruction of the shape of human wart virus
  using spatial correlations.
\newblock {\em Ultramicroscopy}, 17(3):251--262, 1985.

\bibitem{resolution_revolution}
W.~K\"{u}hlbrandt.
\newblock The resolution revolution.
\newblock {\em Science}, 343:1443--1444, 2014.

\bibitem{Kurta2017}
R.~P. Kurta, J.~J. Donatelli, C.~H. Yoon, P.~Berntsen, J.~Bielecki, B.~J.
  Daurer, H.~DeMirci, P.~Fromme, M.~F. Hantke, F.~R. Maia, et~al.
\newblock Correlations in scattered x-ray laser pulses reveal nanoscale
  structural features of viruses.
\newblock {\em Physical review letters}, 119(15):158102, 2017.

\bibitem{Natterer2001a}
F.~Natterer.
\newblock {\em {The Mathematics of Computerized Tomography}}.
\newblock Society for Industrial and Applied Mathematics, 2001.

\bibitem{Punjani2017}
A.~Punjani, J.~L. Rubinstein, D.~J. Fleet, and M.~A. Brubaker.
\newblock {cryoSPARC: algorithms for rapid unsupervised cryo-EM structure
  determination}.
\newblock {\em Nat. Meth.}, 14(3):290--296, 2017.

\bibitem{Radermacher2}
M.~Radermacher, T.~Wagenknecht, A.~Verschoor, and J.~Frank.
\newblock {Three-dimensional structure of the large ribosomal subunit from
  Escherichia coli}.
\newblock {\em EMBO J.}, 6(4):1107--14, 1987.

\bibitem{Saldin2011}
D.~K. Saldin, H.-C. Poon, P.~Schwander, M.~Uddin, and M.~Schmidt.
\newblock {Reconstructing an icosahedral virus from single-particle diffraction
  experiments}.
\newblock {\em Opt. Express}, 19(18):17318, aug 2011.

\bibitem{Saldin2010}
D.~K. Saldin, V.~L. Shneerson, M.~R. Howells, S.~Marchesini, H.~N. Chapman,
  M.~Bogan, D.~Shapiro, R.~A. Kirian, U.~Weierstall, K.~E. Schmidt, et~al.
\newblock {Structure of a single particle from scattering by many particles
  randomly oriented about an axis: toward structure solution without
  crystallization?}
\newblock {\em New J. Phys.}, 12(3):035014, 2010.

\bibitem{Salzman1990}
D.~B. Salzman.
\newblock A method of general moments for orienting {2D} projections of unknown
  {3D} objects.
\newblock {\em Comput. Vision Graph.}, 50:129--156, 1990.

\bibitem{relion}
S.~H.~W. Scheres.
\newblock {RELION:} implementation of a {Bayesian} approach to {cryo-EM}
  structure determination.
\newblock {\em J. Struct. Biol.}, 180(3):519 -- 530, 2012.

\bibitem{Schonemann1966}
P.~H. Schönemann.
\newblock {A generalized solution of the orthogonal procrustes problem}.
\newblock {\em Psychometrika}, 31(1):1--10, 1966.

\bibitem{Vainshtein}
B.~Vainshtein and A.~Goncharov.
\newblock {Determination of the spatial orientation of arbitrarily arranged
  identical particles of unknown structure from their projections}.
\newblock {\em Soviet Physics Doklady}, 31:278, 1986.

\bibitem{VanHeel}
M.~{van Heel}.
\newblock Angular reconstitution: A posteriori assignment of projection
  directions for 3d reconstruction.
\newblock {\em Ultramicroscopy}, 21(2):111 -- 123, 1987.

\end{thebibliography}

\end{document}